# Quillen stratification for the Steenrod algebra

By John H. Palmieri*

### Introduction

Let $A$ be the mod 2 Steenrod algebra. Its cohomology, $H^*(A; \mathbf{F}_2) = \operatorname{Ext}_A^*(\mathbf{F}_2, \mathbf{F}_2)$, is the $E_2$-term of the Adams spectral sequence converging to the 2-component of $\pi_*^s S^0$, the stable homotopy groups of spheres; as such, this cohomology algebra has been studied extensively for over thirty-five years. There have been a number of nice structural results and degree-by-degree computations, many of which have implications for $\pi_*^s S^0$. Information about the ring structure of $H^*(A; \mathbf{F}_2)$, though, has been harder to come by.

Analogously, in the study of the stable homotopy groups of spheres, there were many stem-by-stem calculations before Nishida proved in [14] that every element in $\bigoplus_{i>0} \pi_i^s S^0$ was nilpotent. Since $\pi_0^s S^0 \cong \mathbf{Z}$ and $\pi_i^s S^0 = 0$ when $i < 0$, then Nishida's theorem identifies the ring $\pi_*^s S^0$ modulo nilpotent elements. This result led eventually to the nilpotence theorem of Devinatz, Hopkins, and Smith [3], which has itself led to tremendous structural information about the stable homotopy category.

For another analogy, let $G$ be a compact group and let $k$ be a field of characteristic $p > 0$. In his landmark 1971 paper [18], Quillen gave a description of the cohomology algebra $H^*(G; k)$ modulo nilpotent elements as an inverse limit of the cohomology algebras of the elementary abelian $p$-subgroups of $G$. This has led, for instance, to the work of Benson, Carlson, et al. on the theory of varieties for $kG$-modules, and in general to deep structural information about modular representations of finite groups.

In this paper we describe the cohomology of the Steenrod algebra, modulo nilpotent elements, according to the recipe in Quillen's result: as the inverse limit of the cohomology rings of the "elementary abelian" sub-Hopf algebras of $A$. We hope that this leads to further study in several directions. First of all, according to the point of view of axiomatic stable homotopy theory in [7], the representation theory of any cocommutative Hopf algebra has formal similarities to stable homotopy theory; viewed this way, our main result is an

*Research partially supported by National Science Foundation grant DMS-9407459.

analogue of Nishida's theorem, and one can hope that it will lead to results like the thick subcategory theorem of Hopkins and Smith. For instance, see Conjecture 1.4 for a suggested classification of thick subcategories of finite $A$-modules. Secondly, since our main result is a close analogue of Quillen's theorem, then it provides evidence that other structural results from group theory, such as the study of varieties for modules, may carry over to results about modules over the Steenrod algebra. Thirdly, since the Steenrod algebra is intimately tied to stable homotopy theory via mod 2 cohomology and the Adams spectral sequence, then any new information about $A$ could potentially have topological implications.

The proof is also worth noting. Since $A$ is not "compact" in any sense—$A$ is infinite-dimensional and not profinite, and $H^*(A; \mathbf{F}_2)$ is not Noetherian—one should not expect the group theoretic proof to carry over. Instead, we turn to axiomatic stable homotopy theory. In the appropriate setting, $H^*(A; \mathbf{F}_2)$ plays the role of the stable homotopy groups of spheres. Homotopy theorists have developed many tools, such as Adams spectral sequences, for studying $\pi_*^s S^0$. The axiomatic approach lets us focus these tools on $H^*(A; \mathbf{F}_2)$. In this paper, we use properties of generalized Adams spectral sequences, in the setting of $A$-modules, to prove our main result. We expect that the stable homotopy theoretic approach to Hopf algebra cohomology and representations can be applied fruitfully to other Hopf algebras.

We should point out that one can also study group representation theory from this point of view; in that setting, the role of $\pi_*^s S^0$ is played by $H^*(G; k)$. Hence Nishida's theorem, Quillen's theorem, and our result are all statements about the ring $\pi_*^s S^0$ modulo nilpotence, in various settings. Stable homotopy theoretic ideas have appeared in some recent work in group theory; see [2], for example.

# 1. Results

The mod 2 Steenrod algebra $A$ is the algebra of stable operations on mod 2 cohomology. It is a graded connected cocommutative Hopf algebra over the field $\mathbf{F}_2$; Milnor showed in [13] that its dual $A_*$ is isomorphic, as an algebra, to a polynomial algebra:

$$(*) \qquad A_* \cong \mathbf{F}_2[\xi_1, \xi_2, \xi_3, \ldots],$$

which is graded by setting $|\xi_n| = 2^n - 1$. The product on $A$ is dual to the diagonal map on $A_*$, which is given by

$$\Delta : \xi_n \longmapsto \sum_{i=0}^{n} \xi_{n-i}^{2^i} \otimes \xi_i.$$

(We take $\xi_0$ to be 1.)

*Notation.* Given a commutative ring $k$, a graded augmented $k$-algebra $\Gamma$, and a graded left $\Gamma$-module $M$, we write $H^*(\Gamma; M)$ for $\operatorname{Ext}_\Gamma^{**}(k, M)$. Note that this is a bigraded $k$-module. We write $M^\Gamma$ for $\operatorname{Hom}_\Gamma^*(k, M) \subseteq M$; then $M^\Gamma$ is the set of *invariants* in $M$ under the $\Gamma$-action. It consists of all $x$ in $M$ which support no nontrivial operations by elements of $\Gamma$.

Note that when $\Gamma$ is a (graded) cocommutative Hopf algebra over a field $k$, then $H^*(\Gamma; k)$ is a (graded) commutative $k$-algebra—one can use the diagonal map on $\Gamma$ to induce the product on $H^*(\Gamma; k)$, so cocommutativity of $\Gamma$ implies commutativity of $H^*(\Gamma; k)$. This is the case when $k = \mathbf{F}_2$ and $\Gamma = A$.

Part (a) of this next definition is due to Quillen [18]; part (b) is due to Wilkerson [22].

*Definition* 1.1. Fix a field $k$ of characteristic $p > 0$.

(a) Given a map $f : R \to S$ of graded commutative $k$-algebras, we say that $f$ is an *F-isomorphism* if

   (i) every $x \in \ker f$ is nilpotent, and
   
   (ii) every $y \in \operatorname{coker} f$ is nilpotent, in the sense that for all $y \in S$, there is an $n$ so that $y^{p^n}$ is in the image of $f$.

(b) A Hopf algebra $E$ over $k$ is *elementary* if $E$ is commutative and cocommutative, and has $x^p = 0$ for all $x$ in $IE$ (where $IE = \ker(E \to k)$ is the augmentation ideal of $E$).

*Remark.*

(a) If one can choose the same $n$ for every $y$ in Definition 1.1(a)(ii), then the $F$-isomorphism is said to be *uniform*. We do not expect this to be the case in our results.

(b) When $\operatorname{char}(k) = 2$, then a cocommutative Hopf algebra over $k$ is elementary if and only if it is exterior, as an algebra. The elementary sub-Hopf algebras of the mod 2 Steenrod algebra $A$ were classified by Lin [8, 1.1]; they have also been studied in [22] and [15]. One should think of them as the "elementary abelian subgroups" of $A$.

Given a Hopf algebra $\Gamma$ over a field $k$, for every sub-Hopf algebra $\Lambda$ of $\Gamma$, there is a *restriction map* $H^*(\Gamma; k) \to H^*(\Lambda; k)$. If we let $\mathcal{Q}$ denote the category of elementary sub-Hopf algebras of $\Gamma$, where the morphisms are inclusions, then we can assemble the restriction maps $H^*(\Gamma; k) \to H^*(E; k)$, for $E$ in $\mathcal{Q}$, into

$$\rho : H^*(\Gamma; k) \to \varprojlim_{\mathcal{Q}} H^*(E; k).$$

In the case $k = \mathbf{F}_2$ and $\Gamma = A$, we show in Section 2 that there is an action of $A$ on $\varprojlim_{\mathcal{Q}} H^*(E; \mathbf{F}_2)$, making the inverse limit an algebra over $A$.

Here is our main result.

THEOREM 1.2. *Let $A$ be the mod 2 Steenrod algebra, and let $\mathcal{Q}$ be the category of elementary sub-Hopf algebras of $A$ with morphisms given by inclusions. The map $\rho : H^*(A; \mathbf{F}_2) \to \varprojlim_{\mathcal{Q}} H^*(E; \mathbf{F}_2)$ factors through*

$$\gamma : H^*(A; \mathbf{F}_2) \to \left( \varprojlim_{\mathcal{Q}} H^*(E; \mathbf{F}_2) \right)^A,$$

*and $\gamma$ is an $F$-isomorphism.*

Partial results related to this theorem can be found in [8], [22], [4], and [15]. Lin proves a special case of Theorem 1.2 (actually a version of Theorem 4.1, from which we deduce Theorem 1.2) for certain sub-Hopf algebras of $A$; in his setting, though, some power of every element is invariant, so one does not have to apply $(-)^A$. The other papers all focus on the detection of nilpotence— finding a map $f$ that satisfies condition (i) in Definition 1.1(a)—rather than on the identification of all of the non-nilpotent elements. In particular, the main results of [4] and [15] each imply that $\rho$, and hence $\gamma$, satisfy 1.1(a)(i). We use this result in the proof of Proposition 4.2.

The statement of Theorem 1.2 is strongly influenced by group theory; it is based, of course, on Quillen's theorem [18, 6.2] which identifies the cohomology $H^*(G; k)$ of a compact group $G$, up to $F$-isomorphism, as the inverse limit $\varprojlim H^*(E; k)$. This inverse limit is taken over the category whose objects are elementary abelian $p$-subgroups of $G$ (where $p = \text{char}(k)$), and whose morphisms are generated by inclusions and conjugations by elements of $G$. Using conjugations is technically unpleasant in our setting, so we just use inclusions as our morphisms; this is why we have to consider the invariants under the $A$-action.

We should point out that Theorem 1.2 will not generalize, as stated, to an arbitrary Hopf algebra, because it is known that elementary Hopf algebras do not always detect nilpotence—i.e., the map $\rho : H^*(\Gamma; k) \to \varprojlim_{\mathcal{Q}} H^*(E; k)$ need not satisfy condition (i) of Definition 1.1(a). This is the case, for example, with the odd primary Steenrod algebra; see [22, 6.5]. One has replacements for elementary Hopf algebras, the so-called "quasi-elementary" Hopf algebras, which do a better job of detecting nilpotence. For example, Theorem 1.2 holds with $A$ replaced by any finite-dimensional graded connected cocommutative Hopf algebra $\Gamma$, if one takes the inverse limit over the category of sub-Hopf algebras of quasi-elementary sub-Hopf algebras of $\Gamma$ [16, Theorem 1.4]; one can prove this by imitating a proof of Quillen's theorem for finite groups. See [22] and [16] for a further discussion of these issues.

Here is our second result: we can compute the inverse limit in Theorem 1.2, and we have a formula for the $A$-action.

THEOREM 1.3. *There is an isomorphism of $A$-algebras*

$$\varprojlim_{\mathcal{Q}} H^*(E; \mathbf{F}_2) \cong \mathbf{F}_2[h_{ts} \mid s < t] / (h_{ts}h_{vu} \mid u \geq t),$$

*where $|h_{ts}| = (1, 2^s(2^t - 1))$. The action of $A$ is given by the following on the generators*:

$$\mathrm{Sq}^{2^k}(h_{ts}) = \begin{cases} h_{t-1,s+1} & \text{if } k = s \text{ and } s+1 < t-1, \\ h_{t-1,s} & \text{if } k = s+t-1 \text{ and } s < t-1, \\ 0 & \text{otherwise,} \end{cases}$$

*and then by extending multiplicatively.*

The polynomial generator $h_{ts}$ is represented by $[\xi_t^{2^s}]$ in the cobar construction for any elementary $E$ for which this makes sense. By "one extends multiplicatively," we mean that one uses the Cartan formula to determine the action of $\mathrm{Sq}^n$ on products.

This theorem allows one to compute some invariant elements in the inverse limit, and hence to predict the existence of families of elements in the cohomology of the Steenrod algebra. In Section 2 we present some examples of invariants and describe what is known about their lifts to $H^*(A; \mathbf{F}_2)$.

Aside from this computational information, Theorem 1.2 is potentially interesting from a structural point of view. In particular, it has led us to a similar result for Ext with nontrivial coefficients (Theorem 1.5 below), as well as a proposed classification of the thick subcategories of finite $A$-modules (Conjecture 1.4). To state these, we need a bit of notation.

Let $D$ be the sub-Hopf algebra of $A$ generated by the elementary sub-Hopf algebras of $A$; equivalently, $D$ is generated by the $P_t^s$'s with $s < t$; equivalently, $D$ is dual to the following quotient Hopf algebra of $A_*$:

$$D_* = A_*/(\xi_1^2, \xi_2^4, \xi_3^8, \ldots, \xi_n^{2^n}, \ldots).$$

There is an action of $A$ on $H^*(D; \mathbf{F}_2)$ (see Section 4). In Theorem 4.1, we show that the map

$$\varphi : H^*(A; \mathbf{F}_2) \to H^*(D; \mathbf{F}_2)^A$$

is an $F$-isomorphism; this is an important intermediate result in proving Theorem 1.2.

Recall that if $\mathcal{C}$ is a triangulated category (as defined in [10] or [7], for example), then a subcategory $\mathcal{D} \subseteq \mathcal{C}$ is *thick* (also known as "épaisse") if $\mathcal{D}$ is a triangulated subcategory of $\mathcal{C}$ which is closed under retracts. A classification of the thick subcategories of a given category carries a considerable amount of information; for example, Hopkins and Smith achieved this for the category of finite spectra in [6], and their result has been quite useful in stable homotopy theory. In our situation, the category of $A$-modules is abelian, not triangulated, but neglecting this for a moment, here is our thick subcategory conjecture.

CONJECTURE 1.4. *The thick subcategories of finite $A$-modules are in one-to-one correspondence with the finitely generated radical ideals of $H^*(D; \mathbf{F}_2)$ which are invariant under the action of $A$.*

To get around the abelian/triangulated problem, one should work in a nice triangulated category in which the category of finite $A$-modules sits; then Conjecture 1.4 classifies the thick subcategories of finite objects of that. One good option for such a category is the category of cochain complexes of injective $A_*$-comodules. See [17] for details.

Given a finite $A$-module $M$, we let $I(M)$ be the radical of the annihilator ideal of $H^*(D; M)$ as an $H^*(D; \mathbf{F}_2)$-module, under the Yoneda composition action. The maps in the proposed bijection should go like this: given an invariant radical ideal $I \subseteq H^*(D; \mathbf{F}_2)$, then we let $\mathcal{C}$ be the full subcategory with objects

$$\{M \mid M \text{ finite, } I(M) \supseteq I\}.$$

This is clearly thick. Given a thick subcategory $\mathcal{C}$, then we let $I$ be the ideal

$$\bigcap_{M \in \mathrm{ob}\,\mathcal{C}} I(M).$$

With a little work, one can show that this is an invariant ideal of $H^*(D; \mathbf{F}_2)$ (and we clearly have $I = \sqrt{I}$).

Conjecture 1.4 is supported by the following result, which is proved elsewhere. Part (b) is an analogue of the result of Hopkins and Smith [6, Theorem 11], in which they identify the center of $[X, X]$ up to $F$-isomorphism, for any finite $p$-local spectrum $X$. We write $i$ for the following map:

$$i : \left(\mathrm{Ext}_D^{**}(\mathbf{F}_2, \mathbf{F}_2)/I(M)\right)^A \to \mathrm{Ext}_D^{**}(M, M).$$

THEOREM 1.5. *Let $M$ be a finite-dimensional $A$-module.*

(a) [15] *If $x \in \mathrm{Ext}_A^{**}(M, M)$ is in the kernel of the restriction map*

$$\mathrm{Ext}_A^{**}(M, M) \to \mathrm{Ext}_D^{**}(M, M),$$

*then $x$ is nilpotent.*

(b) [17] *For every*

$$y \in \left(\mathrm{Ext}_D^{**}(\mathbf{F}_2, \mathbf{F}_2)/I(M)\right)^A,$$

*there is an $n$ so that $i(y)^{2^n}$ is in the image of the restriction map*

$$\operatorname{Ext}_A^{**}(M, M) \to \operatorname{Ext}_D^{**}(M, M).$$

*In fact, $i(y)^{2^n}$ is hit by a central element of $\operatorname{Ext}_A^{**}(M, M)$.*

One might in fact conjecture that there is an $F$-isomorphism between the center of $\operatorname{Ext}_A^{**}(M, M)$ and $\left(\operatorname{Ext}_D^{**}(\mathbf{F}_2, \mathbf{F}_2)/I(M)\right)^A$.

Here is the structure of this paper. In Section 2, we discuss elementary and normal sub-Hopf algebras of the Steenrod algebra; we describe the action of $A$ on $\varprojlim H^*(E; \mathbf{F}_2)$, and we prove Theorem 1.3. We also provide some examples of invariant elements in the inverse limit, as well as their lifts to $H^*(A; \mathbf{F}_2)$. In Section 3, we give two results about the spectral sequence associated to a Hopf algebra extension: first, the property of having a vanishing plane with fixed normal vector, with some intercept and at some $E_r$-term, is generic. This is a translation, into algebraic terms, of a stable homotopy theory result: the spectral sequence in question can be viewed as a generalized Adams spectral sequence in the appropriate setting, and this result was originally proved for Adams spectral sequences by Hopkins and Smith; see [5]. The second result in Section 3 is that if we are working with Hopf algebras over the field $\mathbf{F}_2$, for any class $y \in E_2^{0,t}$, then $y^{2^i}$ survives to the $E_{2^i+1}$-term.

In Section 4, we introduce a sub-Hopf algebra $D$ of $A$, and show that Theorem 1.2 follows if we can prove a similar result—Theorem 4.1—about the restriction map $H^*(A; \mathbf{F}_2) \to H^*(D; \mathbf{F}_2)$. Section 5 contains the heart of the proof. Essentially, we examine the Lyndon-Hochschild-Serre spectral sequence associated to the extension $D \to A \to A//D$; we produce a module $M$ for which the $E_2$-term of the spectral sequence converging to $H^*(A; M)$ has a nice vanishing plane, and then we use properties of $M$, along with genericity of vanishing planes and some Ext computations of Lin [8], to deduce that the spectral sequence converging to $H^*(A; \mathbf{F}_2)$ has a similar vanishing plane at the $E_r$-term for some $r$. Theorem 4.1 follows easily from this: if $x \in H^*(A; \mathbf{F}_2)$ is in the kernel of the restriction map to $H^*(D; \mathbf{F}_2)$, then one can show that all sufficiently large powers of $x$ will be above the vanishing plane, and hence zero. For any $y \in (H^*(D; \mathbf{F}_2))^A$, then $y$ represents a class $E_2^{0,t}$ for some $t$, and hence $y^{2^i}$ survives to the $E_{2^i+1}$-term. One can show that for all sufficiently large $i$, the targets of the differentials on $y^{2^i}$ will be above the vanishing plane, and hence zero, so $y^{2^i}$ survives the spectral sequence for $i \gg 0$.

*Remark.* All of the main results of this paper—Theorems 1.2, 1.3, and 4.1—have obvious analogues with $A$ replaced by any sub-Hopf algebra $B$ of $A$, and $\mathcal{Q}$ by the category of elementary sub-Hopf algebras of $B$ (so the objects will be $B \cap E$, for $E \in \operatorname{ob} \mathcal{Q}$). The proofs that we give for the $B = A$ case carry over easily to this slightly more general setting; we omit the details.

*Convention.* Throughout the paper, we have a number of results about the Steenrod algebra, its sub-Hopf algebras, and modules over it. We let $A$ denote the mod 2 Steenrod algebra, and we use letters like $B$, $C$, $D$, and $E$ to denote sub-Hopf algebras of it. We write $B_*$ (etc.) for the corresponding quotient Hopf algebras of the dual of $A$, $A_*$. (We use subscripts for the duals because they are more closely affiliated with homology; for instance, $A_*$ is the mod 2 homology of the mod 2 Eilenberg-Mac Lane spectrum.) Now, we also have a number of results about general graded cocommutative Hopf algebras over a field $k$. To set such results apart, we will use capital Greek letters, like $\Gamma$ and $\Lambda$, to denote such Hopf algebras; typically, $\Lambda$ will denote a sub-Hopf algebra of $\Gamma$. Whether we are working over $A$ or $\Gamma$, "module" means left module; also, all modules are assumed to be graded.

## 2. Sub-Hopf algebras of $A$

In this section we discuss two kinds of sub-Hopf algebras of $A$: elementary ones and normal ones. We also construct the action of $A$ on $\varprojlim H^*(E; \mathbf{F}_2)$, we prove Theorem 1.3, and we use that theorem to find some examples of invariant elements in the inverse limit.

First we bring the reader up to speed on elementary sub-Hopf algebras of $A$. Equation (∗) in Section 1 says that the dual of $A$ is isomorphic to $\mathbf{F}_2[\xi_1, \xi_2, \ldots]$; we dualize with respect to the monomial basis of $A_*$, and we let $P_t^s$ be the basis element dual to $\xi_t^{2^s}$.

PROPOSITION 2.1.

(a) *As an algebra, each elementary sub-Hopf algebra $E$ of $A$ is isomorphic to the exterior algebra on the $P_t^s$'s that it contains.*

(b) *Hence,*
$$H^*(E; \mathbf{F}_2) \cong \mathbf{F}_2[h_{ts} \mid P_t^s \in E],$$
*where $h_{ts}$ is represented by $[\xi_t^{2^s}]$ in the cobar complex for $E$; in particular, the degree of $h_{ts}$ is $|h_{ts}| = (1, 2^s(2^t - 1))$.*

(c) *The maximal elementary sub-Hopf algebras of $A$ are of the form*
$$\left( \mathbf{F}_2[\xi_1, \xi_2, \ldots]/(\xi_1, \ldots, \xi_{n-1}, \xi_n^{2^n}, \xi_{n+1}^{2^n}, \xi_{n+2}^{2^n}, \ldots) \right)^*,$$
*for $n \geq 1$.*

*Proof.* This is an accumulation of well-known facts about sub-Hopf algebras of the Steenrod algebra. Part (b) follows immediately from part (a); parts (a) and (c) can be found in [10, pp. 233–235]. □

We also need to discuss normal sub-Hopf algebras.

*Definition* 2.2.  Suppose that $\Gamma$ is a cocommutative Hopf algebra over a field $k$. A sub-Hopf algebra $\Lambda$ of $\Gamma$ is *normal* if the left ideal of $\Gamma$ generated by $\Lambda$ is equal to the right ideal of $\Gamma$ generated by $\Lambda$ (i.e., we have $\Gamma \cdot I\Lambda = I\Lambda \cdot \Gamma$, where $I\Lambda = \ker(\Lambda \to k)$ is the augmentation ideal of $\Lambda$). If $\Lambda$ is normal, then we have $\Gamma \otimes_\Lambda k = k \otimes_\Lambda \Gamma$, and this vector space has the structure of a Hopf algebra; we use $\Gamma//\Lambda$ to denote it. In this situation, we say that $\Lambda \to \Gamma \to \Gamma//\Lambda$ is a *Hopf algebra extension*.

First we discuss properties of normal sub-Hopf algebras of an arbitrary graded cocommutative Hopf algebra $\Gamma$, and then we specialize to the case of normal sub-Hopf algebras of the Steenrod algebra $A$.

PROPOSITION 2.3.  *Let $\Gamma$ be a graded cocommutative Hopf algebra over a field $k$. Fix a $\Gamma$-module $M$ and a normal sub-Hopf algebra $\Lambda$ of $\Gamma$.*

(a) *For each integer $t$, there is an action of $\Gamma//\Lambda$ on $H^{t,*}(\Lambda; M)$. This action is natural in $M$, and if viewed as an action of $\Gamma$ on $H^{t,*}(\Lambda; M)$, then it is natural in $\Lambda$ as well. This action lowers internal degrees: given $\gamma \in \Gamma//\Lambda$ of degree $i$ and $x \in H^{t,v}(\Lambda; M)$, then $\gamma x$ is in $H^{t,v-i}(\Lambda; M)$.*

(b) *There is a spectral sequence, the* Lyndon-Hochschild-Serre *spectral sequence, with*

$$E_2^{s,t,u}(M) = H^{s,u}(\Gamma//\Lambda; H^{t,*}(\Lambda; M)) \Rightarrow H^{s+t,u}(\Gamma; M),$$

*and with differentials*

$$d_r : E_r^{s,t,u} \to E_r^{s+r,t-r+1,u}.$$

*It is natural in $M$.*

*Proof.* This is a simple extension of the corresponding results for group extensions. For Hopf algebra extensions, this spectral sequence is discussed in [20] and many other places. (We also outline a construction of the spectral sequence in the proof of Theorem 3.2.) □

By the way, we will often be applying Proposition 2.3(a) in the case where $\Gamma = A$ is the mod 2 Steenrod algebra, yielding an action of $A$ on $H^*(B; \mathbf{F}_2)$ for $B$ a normal sub-Hopf algebra of $A$. Note that one also has Steenrod operations acting on the cohomology of any cocommutative Hopf algebra over $\mathbf{F}_2$; see May's paper [12], for example. These two actions are *quite different*; in this paper, we will almost always use the $A$-action arising as in Proposition 2.3(a) (the only exception is in the proof of Lemma 3.7).

PROPOSITION 2.4.  *Let A be the mod* 2 *Steenrod algebra.*

(a) *A sub-Hopf algebra of A is normal if and only if its dual is of the form*

$$\mathbf{F}_2[\xi_1, \xi_2, \ldots]/(\xi_1^{2^{n_1}}, \xi_2^{2^{n_2}}, \xi_3^{2^{n_3}}, \ldots),$$

*with* $n_1 \leq n_2 \leq n_3 \leq \ldots \leq \infty$. (*If some* $n_i = \infty$ *in this expression, this means that one does not divide out by any power of* $\xi_i$.)

(b) *Hence the maximal elementary sub-Hopf algebras of A are normal.*

*Proof.* This is well-known; see [10, Theorem 15.6], for example.  □

Now we construct the action of $A$ on $\varprojlim H^*(E; \mathbf{F}_2)$. We let $\mathcal{Q}_0$ denote the full subcategory of $\mathcal{Q}$ consisting of the normal elementary sub-Hopf algebras of $A$. By Proposition 2.4(b), we see that $\mathcal{Q}_0$ is final in $\mathcal{Q}$; hence we have

$$\varprojlim_{\mathcal{Q}} H^*(E; \mathbf{F}_2) = \varprojlim_{\mathcal{Q}_0} H^*(E; \mathbf{F}_2).$$

On the right we have an inverse limit of modules over $A$; we give this the induced $A$-module structure. So we have (since taking invariants is an inverse limit)

$$\left(\varprojlim_{\mathcal{Q}} H^*(E; \mathbf{F}_2)\right)^A = \left(\varprojlim_{\mathcal{Q}_0} H^*(E; \mathbf{F}_2)\right)^A$$
$$= \varprojlim_{\mathcal{Q}_0} \left(H^*(E; \mathbf{F}_2)^A\right)$$
$$= \varprojlim_{\mathcal{Q}_0} \left(H^*(E; \mathbf{F}_2)^{A/\!/E}\right).$$

Theorem 1.3 follows more or less immediately:

*Proof of Theorem* 1.3. We want to compute $\varprojlim H^*(E; \mathbf{F}_2)$, and we want to compute the action of $A$ on this ring. Proposition 2.1(b) tells us $H^*(E; \mathbf{F}_2)$ for each $E$; since the maps in the category $\mathcal{Q}$ are inclusion maps of one exterior algebra into another exterior algebra in which one simply adds more generators, the induced maps on cohomology are surjections of one polynomial ring onto another, each polynomial generator mapping either to zero or to a polynomial generator (of the same name). So assembling the cohomology rings of the maximal elementary sub-Hopf algebras of $A$ gives our computation of the inverse limit.

To get our formula for the action of $A$ on $\varprojlim_{\mathcal{Q}} H^*(E; \mathbf{F}_2)$, we have to compute the $A$-action on $H^*(E; \mathbf{F}_2)$ for $E \leq A$ maximal elementary. This is done in [15, Lemmas A.3 and A.5]. (Briefly: for degree reasons, we have

$$\mathrm{Sq}^{2^k}(h_{ts}) = \begin{cases} \alpha h_{t-1,s+1} & \text{if } k = s \text{ and } s+1 < t-1, \\ \beta h_{t-1,s} & \text{if } k = s+t-1 \text{ and } s < t-1, \\ 0 & \text{otherwise}, \end{cases}$$

for some $\alpha, \beta \in \mathbf{F}_2$. So one only has to check that $\alpha = \beta = 1$; this is a reasonably straightforward computation.) □

One can rewrite the action of $A$ on the inverse limit as a coaction of the dual $A_*$—that is, as a map $\varprojlim H^*(E; \mathbf{F}_2) \to A_* \otimes \varprojlim H^*(E; \mathbf{F}_2)$. The coaction is multiplicative, and its value on the polynomial generators is given by the following formula:

$$h_{ts} \longmapsto \sum_{j=0}^{\lfloor \frac{s+t-1}{2} \rfloor} \sum_{i=j+s+1}^{t-j} \zeta_j^{2^s} \xi_{t-i-j}^{2^{i+j+s}} \otimes h_{i,j+s}.$$

(As usual, $\zeta_j$ denotes $\chi(\xi_j)$, where $\chi : A_* \to A_*$ is the conjugation map.) To verify this formula, one checks that it is correct on the primitives in $A_*$—the terms of the form $\xi_1^{2^k} \otimes h_{m,n}$—and that it is coassociative. Milnor computed $\chi(\xi_j)$ as a polynomial in the $\xi_i$'s in [13, Lemma 10]; in particular, he found that

$$\chi(\xi_j) = \xi_1^{2^j - 1} + \text{(other terms)}.$$

One can use this observation to get a formula for the action of every $\mathrm{Sq}^n$ in $A$, not just the indecomposables:

$$\mathrm{Sq}^n(h_{ts}) = \begin{cases} h_{t-j,s+j} & \text{if } n = 2^{s+j} - 2^s \text{ and } s+j < t-j, \\ h_{t-j-1,s+j} & \text{if } n = 2^{s+t-1} + 2^{s+j} - 2^s \text{ and } s+j < t-j-1, \\ 0 & \text{otherwise}. \end{cases}$$

Using this explicit description of the (co)action, one can find elements in the ring of invariants, either by hand or by computer; we give a few examples below. We have not been able to identify the entire ring of invariants, though. For these examples, we let $R$ denote the bigraded ring $\mathbf{F}_2[h_{ts} \mid s < t]/(h_{ts}h_{vu} \mid u \geq t)$. We give a graphical depiction of the (co)action of $A$ on $R$ in Figure 1.

*Example.*
(a) The element $h_{t,t-1}$ in $R^{1,2^{t-1}(2^t-1)}$ is an invariant, for $t \geq 1$. Indeed, we know that $h_{10}$ lifts to $H^{1,1}(A; \mathbf{F}_2)$; also $h_{21}^4$ lifts to an element in $H^{4,24}(A; \mathbf{F}_2)$ (the element known as $g$ or $\overline{\kappa}$; see [23] for the identification of this element as a lift of a power of $h_{21}$). We do not know which power of $h_{t,t-1}$ survives for $t \geq 3$.

(b) We have some families of invariants. Note that $h_{20}$ is not invariant—it supports a single nontrivial operation:

$$\mathrm{Sq}^2(h_{20}) = h_{10}.$$

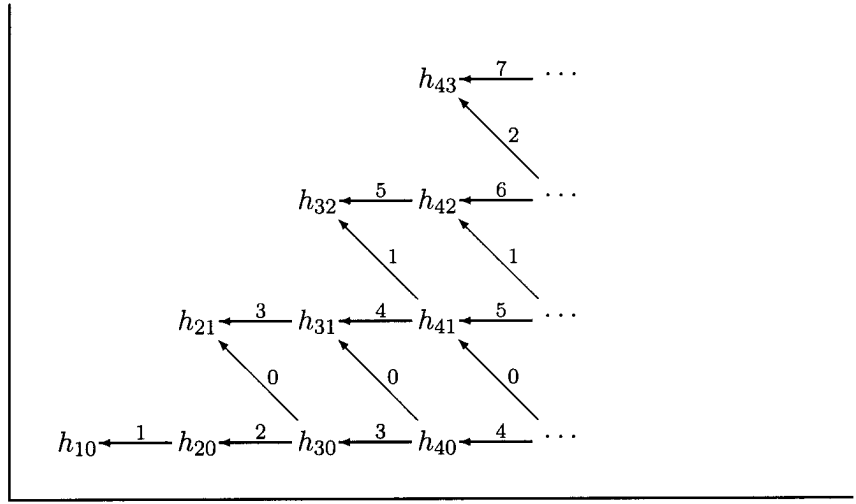

**Figure 1.** Graphical depiction of the action of $A$ on $R = \mathbf{F}_2[h_{ts} \mid s < t] / (h_{ts}h_{vu} \mid u \geq t)$. An arrow labeled by a number $k$ represents an action by $\text{Sq}^{2^k} \in A$, or equivalently a "coaction" by $\xi_1^{2^k} \in A_*$—in other words, a term of the form $\xi_1^{2^k} \otimes$ (target) in the coaction on the source.

Then since $h_{21}$ is invariant, the element $h_{20}h_{21}$ supports at most one possibly nontrivial operation:

$$\text{Sq}^2(h_{20}h_{21}) = h_{10}h_{21}.$$

But since $h_{10}h_{21} = 0$ in $R$, then $h_{20}h_{21}$ is invariant; indeed, $h_{20}^i h_{21}^j$ is invariant for all $i \geq 0$ and $j \geq 1$. It turns out that more of these lift to $H^*(A; \mathbf{F}_2)$ than one might expect from Theorem 1.2:

| element of $R$ | bidegree | lifts to | element of $R$ | bidegree | lifts to |
|---|---|---|---|---|---|
| $h_{20}^2 h_{21}^2$ | (4,18) | $d_0$ | $h_{20}^4 h_{21}^3$ | (7,30) | $i$ |
| $h_{20} h_{21}^3$ | (4,21) | $e_0$ | $h_{20}^3 h_{21}^4$ | (7,33) | $j$ |
| $h_{21}^4$ | (4,24) | $g$ | $h_{20}^2 h_{21}^5$ | (7,36) | $k$ |

(See [23] for the analysis of the elements $d_0$, $e_0$, and $g$ of $H^*(A; \mathbf{F}_2)$; see [11, p. 46] for the elements $i$, $j$, and $k$.)

Furthermore, the "Mahowald-Tangora wedge" [9] contains lifts of the elements $h_{20}^i h_{21}^{8j}$ for all $i \geq 0$ and $j \geq 1$. (See [11]; earlier, Zachariou [23] verified this for elements of the form $h_{20}^{2i} h_{21}^{2(i+j)}$ for $i, j \geq 0$.) These elements are distributed over a wedge between lines of slope $\frac{1}{2}$ and $\frac{1}{5}$ (in the usual Adams spectral sequence picture of $H^*(A; \mathbf{F}_2)$).

(c) Similarly, while $h_{30}$ and $h_{31}$ are not invariant, the monomials

$$\{h_{30}^i h_{31}^j h_{32}^k \mid i \geq 0, \ j \geq 0, \ k \geq 1\}$$

are invariant elements. Hence some powers of them lift to $H^*(A; \mathbf{F}_2)$; we do not know which powers, though. Continuing in this pattern, we find that for $n \geq 1$, we have sets of invariant elements

$$\{h_{n0}^{i_0} h_{n1}^{i_1} \ldots h_{n,n-1}^{i_{n-1}} \,|\, i_0, \ldots, i_{n-2} \geq 0,\ i_{n-1} \geq 1\}.$$

We do not know which powers of these elements lift. In any case, we can see that the family of elements in the Mahowald-Tangora wedge is not a unique phenomenon—we have infinitely many such families, and all but the $n = 1$ and $n = 2$ families give more than a lattice of points in the bigraded vector space $H^*(A; \mathbf{F}_2)$. One can easily show, by the way, that any invariant monomial in $R$ is of the form $h_{n0}^{i_0} h_{n1}^{i_1} \ldots h_{n,n-1}^{i_{n-1}}$ for some $n$. The invariant *polynomials* may be much more complicated.

(d) Margolis, Priddy, and Tangora [11, p. 46] have found some other non-nilpotent elements, such as $x' \in H^{10,63}(A; \mathbf{F}_2)$ and $B_{21} \in H^{10,69}(A; \mathbf{F}_2)$. These both come from the invariant

$$z = h_{40}^2 h_{21}^3 + h_{20}^2 h_{21}^2 h_{41} + h_{30}^2 h_{21} h_{31}^2 + h_{20}^2 h_{31}^3$$

in $R$. While we do not know which power of $z$ lifts to $H^*(A; \mathbf{F}_2)$, we do know that $B_{21}$ maps to the product $h_{20}^3 h_{21}^2 z$, and $x'$ maps to $h_{20}^5 z$. In other words, we have at least found some elements in $(z) \cap R^A$ which lift to $H^*(A; \mathbf{F}_2)$.

(e) Computer calculations have led us to a few other such "sporadic" invariant elements (i.e., invariant elements that do not belong to any family—any family of which we are aware, anyway):

$$h_{20}^8 h_{31}^4 + h_{30}^8 h_{21}^4 + h_{21}^{11} h_{31}$$

in $R^{12,80}$, an element in $R^{13,104}$ (a sum of 12 monomials in the variables $h_{i,0}$ and $h_{i,1}$, $2 \leq i \leq 5$), and an element in $R^{9,104}$ (a sum of eight monomials in the same variables). We do not know what powers of these elements lift to $H^*(A; \mathbf{F}_2)$, nor are we aware of any elements in the ideal that they generate which are in the image of the restriction map from $H^*(A; \mathbf{F}_2)$.

## 3. Extensions of Hopf algebras

In this section we discuss normal sub-Hopf algebras and the Lyndon-Hochschild-Serre spectral sequence that arises from a Hopf algebra extension. In particular, we show in Theorem 3.2 that the presence of a vanishing plane at some term of such a spectral sequence is a "generic" property (Definition 3.4). We also point out in Lemma 3.7 that if $y$ is in the $E_2^{0,t,u}$-term of this spectral sequence, then $y^{2^i}$ survives to the $E_{2^i+1}$-term (if we're working over the field $\mathbf{F}_2$).

(Recall that we defined the notions of normal sub-Hopf algebra and Hopf algebra extension in Definition 2.2, and we presented the spectral sequence of a Hopf algebra extension in Proposition 2.3.)

Theorem 3.2 is our main theorem for this section. It is essentially due to Hopkins and Smith and can be found in [5], in the setting of the Adams spectral sequence associated to a generalized homology theory. Since the Lyndon-Hochschild-Serre spectral sequence can be viewed as an Adams spectral sequence in the appropriate context, it is not surprising that their result applies here.

*Remark.* The "appropriate context" is the category whose objects are unbounded chain complexes of injective comodules over the dual of the Hopf algebra $\Gamma$, and whose morphisms are chain homotopy classes of maps—this is a stable homotopy category [7, 9.5]. We have stated our results and proofs in the language of modules, rather than chain complexes, but if one is willing to work in the chain complex setting, then some of the following can be stated and proved a bit more cleanly (see [17]).

We are going to study properties of the Lyndon-Hochschild-Serre spectral sequence
$$H^{s,u}(\Gamma//\Lambda; H^{t,*}(\Lambda; M)) \Rightarrow H^{s+t,u}(\Gamma; M)$$
by looking at an injective resolution of the module $M$. This is not the usual procedure when working over the Steenrod algebra—often one restricts to bounded-below $A$-modules, in which case there are not enough injectives. To ensure that injective modules are well-behaved (e.g., $I \otimes M$ should be injective when $I$ is), we make some restrictions.

*Convention.* For the remainder of this section, we work with a rather restricted class of modules. First, we work with finite-type modules: modules which are finite-dimensional in each degree. Second, all modules are assumed to be bounded-above—i.e., they are zero in all sufficiently large gradings. We will occasionally refer to bounded-above finite-type modules as "nice," although one should assume that all modules in this section satisfy these restrictions, whether labeled as nice or not. We also assume that the Hopf algebra $\Gamma$ is itself finite-type, so that its dual $\Gamma_*$ is nice.

(Alternatively, one can work with "locally finite" $\Gamma$-modules: modules $M$ so that for every $m \in M$, the submodule generated by $m$ is finite dimensional. The category of locally finite modules is isomorphic to the category of $\Gamma_*$-comodules, and injective comodules are well-behaved.)

Given two $\Gamma$-modules $M$ and $N$, $M \otimes N$ indicates the $\Gamma$-module with the usual diagonal $\Gamma$-action. Similarly, $\text{Hom}_k(M, N)$ is given a $\Gamma$-module structure using conjugation; in particular, the vector space dual of a $\Gamma$-module can be given the structure of a $\Gamma$-module.

LEMMA 3.1. *Suppose that $\Gamma$ is a graded connected cocommutative finite-type Hopf algebra over a field $k$.*

(a) *The dual $\Gamma_*$ of $\Gamma$ is an injective module.*

(b) *If $M$ is a bounded-above finite-type module, and if $I$ is a bounded-above finite-type injective module, then $I \otimes M$ is a bounded-above finite-type injective.*

*Proof.* Part (a): By [10, 11.12], $\Gamma_*$ is injective as a right $\Gamma$-module. The anti-automorphism $\chi$ induces an exact covariant equivalence of categories (with $\chi^2 = 1$), so $\chi$ takes injectives to injectives.

Part (b): An object is injective if and only if mapping into it is exact; thus, if $I$ is injective, so is $\text{Hom}_k(N, I)$ for any $N$. Given the hypotheses on $M$ and $I$, there is an isomorphism $M \otimes I \cong \text{Hom}_k(M_*, I)$, where $M_*$ is the dual of $M$. □

Here is the main result of this section.

THEOREM 3.2. *Suppose that $\Lambda \to \Gamma \to \Gamma//\Lambda$ is an extension of graded connected cocommutative finite-type Hopf algebras over a field $k$. For any bounded-above finite-type $\Gamma$-module $M$, consider the associated spectral sequence*

$$E_2^{s,t,u}(M) = H^{s,u}(\Gamma//\Lambda; H^{t,*}(\Lambda; M)) \Rightarrow H^{s+t,u}(\Gamma; M).$$

*Fix integers $a$, $b$, and $c$, with $a \geq b \geq 0$ and $c < 0$, and consider the following property $P$ of the $\Gamma$-module $M$:*

P: There are numbers $r$ and $d$ so that $E_r^{s,t,u}(M) = 0$ when $as + bt + cu > d$. Then:

(a) *$M$ satisfies $P$ if and only if $\Sigma M$ satisfies $P$.*

(b) *Given a short exact sequence of bounded-above finite-type $\Gamma$-modules*

$$0 \to M' \to M \to M'' \to 0,$$

*if two of $M'$, $M$, and $M''$ satisfy property $P$, then so does the third.*

(c) *If the bounded-above finite-type $\Gamma$-module $L \oplus M$ satisfies property $P$, then so do $L$ and $M$.*

(d) *Every bounded-above finite-type injective $\Gamma$-module satisfies property $P$.*

Here, $\Sigma$ denotes the usual "shift" functor on the module category.

Property $P$ says that for some $r$ and $d$, the $E_r$-term of the spectral sequence has a vanishing plane with normal vector $(a, b, c)$ and intercept $d$. Note that $a$, $b$, and $c$ are fixed throughout, while $r$ and $d$ may vary; for example, the intercept for the module $M$ will be different from the intercept for $\Sigma M$.

*Remark.* We should point out that the result ought to hold without assuming the inequalities on $a$, $b$, and $c$. Indeed, the corresponding result in the chain complex setting of [17] holds without such a restriction, but the proof breaks down into cases, depending on the values of $a$, $b$, and $c$. For our purposes, it suffices to consider the situation as stated in the theorem, so we omit the other cases. See [17, Appendix A.2] for the more general version, stated and proved using homotopy-theoretic language.

*Proof.* Part (d) is clear: if $M$ is a bounded-above injective $\Gamma$-module, then the spectral sequence collapses, and $E_2^{s,t,u} = E_\infty^{s,t,u}$ is zero unless $s = t = 0$ and $u \geq i$, where $i$ is the largest nonzero grading of $M$. Since $c < 0$, then $E_2^{s,t,u}$ is certainly zero above some plane with normal vector $(a, b, c)$.

Parts (a) and (c) are also clear: for part (a), the spectral sequence for $\Sigma M$ is the same as the spectral sequence for $M$, but with a shift in dimension; so $M$ and $\Sigma M$ would have a vanishing plane with the same normal vector at the same $E_r$, with slightly different intercepts. For part (c), the naturality of the spectral sequence ensures that the $E_r$-term for the module $L$ is a retract of the $E_r$-term for $L \oplus M$. So a vanishing plane for the sum induces one on each summand.

We are left with (b). We start by giving a construction of the spectral sequence; this is essentially the same construction as in [1, 3.5], for instance, except we use injective rather than projective resolutions.

Let $J_*$ be a fixed injective resolution of $k$ as a $\Gamma$-module; since $\Gamma$ is assumed to be connected and finite-type, then we may choose this resolution so that each $J_t$ is bounded-above and finite-type. Similarly, let $I_*$ be a fixed resolution of $k$ by nice $\Gamma//\Lambda$-injectives. We view each $I_s$ as a $\Gamma$-module via the quotient map $\Gamma \to \Gamma//\Lambda$. Then we let

$$\begin{aligned} E_0^{s,t,u}(M) &= \mathrm{Hom}_\Gamma^u(k, I_s \otimes J_t \otimes M) \\ &\cong \mathrm{Hom}_{\Gamma//\Lambda}^u(k, \mathrm{Hom}_\Lambda^*(k, I_s \otimes J_t \otimes M)) \\ &\cong \mathrm{Hom}_{\Gamma//\Lambda}^u(k, I_s \otimes \mathrm{Hom}_\Lambda^*(k, J_t \otimes M)). \end{aligned}$$

The first isomorphism holds because of the following lemma. The second holds since $I_s$ is a trivial $\Lambda$-module.

LEMMA 3.3. *For any $\Gamma$-module $N$, there is an isomorphism*

$$\mathrm{Hom}_\Gamma^*(k, N) \cong \mathrm{Hom}_{\Gamma//\Lambda}^*(k, \mathrm{Hom}_\Lambda^*(k, N)).$$

*In other words, $N^\Gamma = (N^\Lambda)^{\Gamma//\Lambda}$.*

*Proof.* This is a straightforward verification. Both of the terms in question are subsets of $N$, and one can easily check that they are equal. □

The maps in the complexes $I_*$ and $J_*$ induce differentials on $E_0^{s,t,u}$, making it a double complex; hence we get two spectral sequences. Since $I_s \otimes J_t \otimes M$ is injective as a $\Gamma$-module by Lemma 3.1, then one of the spectral sequences collapses to $E_2^* = H^*(\Gamma; M)$. Similarly, one can see that the other has $E_2^{s,t,u}$ as in the statement of the theorem.

We need to examine the exact couple associated with this spectral sequence: we let
$$D_0^{s,t,u}(M) = \bigoplus E_0^{i,j,u}(M)$$
where the direct sum is over all $i$ and $j$ with $i + j = s + t$ and $i \geq s$. Then we have

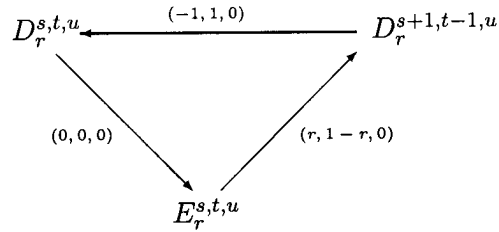

where the arrows are labeled by their degrees (i.e., how much they change the superscripts). The map $i_0 : D_0^{s+1,t-1,u} \to D_0^{s,t,u}$ is the obvious inclusion, and $i_r : D_r^{s+1,t-1,u} \to D_r^{s,t,u}$ is induced by $i_0$, inductively: given $D_{r-1}$ and $i_{r-1}$, then $D_r$ is the image of $i_{r-1}$, and $i_r$ is the restriction of $i_{r-1}$ to this image. The exact couple unfolds into an exact sequence:
$$\ldots \to D_r^{stu} \to E_r^{stu} \to D_r^{s+r,t-r+1,u} \to D_r^{s+r-1,t-r+2,u} \to E_r^{s+r-1,t-r+2,u} \to \ldots.$$

Now, consider the following property:

$Q$: There are numbers $r$ and $d$ so that whenever $as + bt + cu > d$, the composite
$$\underbrace{i_0 \circ \ldots \circ i_0}_{r} : D_0^{s+r,t-r,u} \longrightarrow D_0^{s,t,u}$$
is zero, or equivalently, $D_r^{s,t,u} = 0$.

We claim that properties $P$ and $Q$ are essentially equivalent. We use subscripts, as in $P_{r,d}$, to indicate the values of $r$ and $d$ for which the properties hold. We have the implications
$$P_{r,d} \Rightarrow Q_{r,d+(a-b)(r-1)+b},$$
$$Q_{r,d} \Rightarrow P_{r,d}.$$

These are easy to see: if $M$ satisfies $P$ at the $E_r$-term, with intercept $d$, then we have $E_r^{stu} = 0$ when $as + bt + cu > d$. Since $a \geq b \geq 0$, then we also have $E_r^{s+r-1,t-r+2,u} = 0$ under the same conditions, so the map $i_r : D_r^{s+r,t-r+1,u} \to D_r^{s+r-1,t-r+2}$ is an isomorphism under these conditions. Indeed, one can see that for any integer $k > 0$, the following composite is an isomorphism:

$$(i_r)^k : D_r^{s+r-1+k,t-r+2-k,u} \to D_r^{s+r-1,t-r+2,u}.$$

By construction of the spectral sequence, though, $D_0^{s+n,t-n,u}$ is zero for all $n > t$; since $D_r^{s+n,t-n,u} \subseteq D_0^{s+n,t-n,u}$, the same is true for $D_r^{s+n,t-n,u}$. We conclude that $D_r^{s+r-1,t-r+2,u} = 0$. By the definition of $D_r$, this group is equal to the image of

$$(i_0)^r : D_0^{s+2r-1,t-2r+2,u} \to D_0^{s+r-1,t-r+2,u}.$$

Hence $M$ satisfies property $Q$ at the $E_r$-term, with intercept $d+(a-b)(r-1)+b$.

It is easier to show the other implication: given property $Q$ at $E_r$ with intercept $d$, then when $as + bt + cu > d$, we have $D_r^{stu} = 0 = D_r^{s+r,t-r+1,u}$, so $E_r^{stu} = 0$, by the exact sequence.

By our construction of the spectral sequence, $D_0(-)$ is functorial; furthermore, it takes a short exact sequence of modules to a short exact sequence of double complexes. One can easily see that if two of the three complexes in a short exact sequence satisfy condition $Q$, then so does the third; one may take the intercept $d$ for the third to be the larger of the other two intercepts, and one may take the term $r$ for the third to be the sum of the $r$'s for the other complexes. □

We borrow the following terminology from the analogous situation in a triangulated category (see [6], for example).

*Definition* 3.4. Any property $P$ of bounded-above finite-type $\Gamma$-modules which satisfies (a)–(d) of Theorem 3.2 is called a *generic* property.

*Example* 3.5. Suppose that $\Gamma$ is a cocommutative finite-type Hopf algebra over a field $k$, and let $\Gamma_*$ be the dual of $\Gamma$. If $M$ is a nice $\Gamma$-module, then $\Gamma_* \otimes M$ is a bounded-above finite-type injective, by Lemma 3.1. Hence if $M$ satisfies a generic property $P$, then so does the module

$$\operatorname{coker}(M \hookrightarrow \Gamma_* \otimes M).$$

We let $\Omega^{-1}M$ denote this module. It is standard (and easy to show) that

$$\operatorname{Ext}_\Gamma^{s+1}(\Omega^{-1}L, M) = \operatorname{Ext}_\Gamma^s(L, M),$$
$$\operatorname{Ext}_\Gamma^s(L, \Omega^{-1}M) = \operatorname{Ext}_\Gamma^{s+1}(L, M),$$

for $s > 0$. See [10, 14.8], for example. We point out that if

$$0 \to M' \to M \to M'' \to 0$$

is a short exact sequence of nice modules, then for some injective modules $I$ and $J$, so is the following:

$$0 \to M \to I \oplus M'' \to J \oplus \Omega^{-1} M' \to 0.$$

(See [10, p. 210].) We say that

$$0 \to M \to M'' \to \Omega^{-1} M' \to 0$$

is exact *up to injective summands*.

We will use Theorem 3.2 in combination with the next result.

PROPOSITION 3.6. *Fix a graded connected cocommutative finite-type Hopf algebra $\Gamma$ over a field $k$, and fix a property $P$ which is generic on $\Gamma$-modules. Let*

$$0 \to \Sigma^j N \to M \to N \to 0$$

*be a short exact sequence of bounded-above finite-type $\Gamma$-modules, representing the element*

$$h \in \operatorname{Ext}_\Gamma^{1,j}(N, N).$$

*If $h$ is nilpotent in $\operatorname{Ext}_\Gamma^{**}(N, N)$, then $M$ satisfies $P$ if and only if $N$ satisfies $P$.*

*Proof.* By the definition of genericity, if $N$ satisfies $P$, then so does $M$.

For the converse, we assume that $h^m = 0$. For each integer $i \geq 1$, we define a $\Gamma$-module $M_i$ by the following sequence, which represents $h^i$ in $\operatorname{Ext}_\Gamma^0(N, \Omega^{-i} N) = \operatorname{Ext}_\Gamma^i(N, N)$ and is exact up to injective summands:

$$0 \to M_i \to N \xrightarrow{h^i} \Omega^{-i} N \to 0.$$

(We omit suspensions here—we should write $h_i : N \to \Sigma^{ij} \Omega^{-i} N$.) So $M_i$ is only well-defined up to injective summands. The following is also exact up to injective summands:

$$0 \to \Omega^{-(i-1)} N \to M_i \to N \to 0.$$

Since $h^m = 0$, then this sequence splits when $i = m$; hence we see that $N$ is a summand of $M_m \oplus I$ for some bounded-above injective $I$.

We claim that if $M$ satisfies $P$, then so does $M_i$ for all $i$; verifying this claim (and letting $i = m$) will finish the proof. We prove this by induction on $i$. We start the induction by noting that $M = M_1$.

So suppose that $M_{i-1}$ satisfies $P$. Consider the following commutative diagram:

$$\begin{array}{ccccccccc}
& & 0 & & 0 & & 0 & & \\
& & \downarrow & & \downarrow & & \downarrow & & \\
0 & \to & K & \longrightarrow & 0 & \longrightarrow & \Omega^{-(i-1)}M_1 & \to & 0 \\
& & \downarrow & & \downarrow & \overset{h^{i-1}}{} & \downarrow & & \\
0 & \to & M_{i-1} & \longrightarrow & N & \xrightarrow{h^{i-1}} & \Omega^{-(i-1)}N & \to & 0 \\
& & \downarrow & & \| & & \downarrow \Omega^{-(i-1)}h & & \\
0 & \to & M_i & \longrightarrow & N & \xrightarrow{h^i} & \Omega^{-i}N & \to & 0. \\
& & \downarrow & & \downarrow & & \downarrow & & \\
& & 0 & & 0 & & 0 & &
\end{array}$$

The rows and columns of this diagram are exact up to injective summands, so we see that $K$ is isomorphic to $\Omega^{-(i-2)}M_1$, up to injective summands. By applying Example 3.5, we see that $\Omega^{-(i-2)}M_1$ satisfies $P$, so from the left-hand column we conclude that $M_i$ does as well.  □

We will also need the following result, about the survival of powers of elements in our spectral sequence. Note that when using trivial coefficients, the Lyndon-Hochschild-Serre spectral sequence is a spectral sequence of commutative algebras: each $E_r$-term is a (graded) commutative algebra, and each differential $d_r$ is a derivation. So if we are working with a field $k$ of characteristic $p > 0$, then for any $z \in E_r^{*,*}$, $z^p$ survives to the $E_{r+1}$-term. It turns out that we can greatly improve on this result.

LEMMA 3.7. *Suppose that there is an extension $\Lambda \to \Gamma \to \Gamma//\Lambda$ of graded cocommutative Hopf algebras over the field $\mathbf{F}_2$, and consider the associated spectral sequence $E_r^{s,t,u}(\mathbf{F}_2)$. Fix $y \in E_2^{0,t,u}$. Then for each $i$, the element $y^{2^i}$ survives to $E_{2^i+1}^{0,2^i t, 2^i u}$.*

*Proof.* This follows from properties of Steenrod operations on this spectral sequence, as discussed by Singer in [20]. Suppose that we have a spectral sequence $E_r^{s,t}$ which is a spectral sequence of algebras over the Steenrod algebra. Fix $z \in E_r^{s,t}$, and fix an integer $k$. Then [20, 1.4] tells us to which term of the spectral sequence $\mathrm{Sq}^k z$ survives (the result depends on $r$, $s$, $t$, and $k$). For instance, if $z \in E_r^{0,t}$, then $\mathrm{Sq}^t(z) = z^2$ survives to $E_{2r-1}^{0,2t}$.  □

(There are similar results at odd primes; see [19].)

By the way, Singer's results in [20] are stated in the case of an extension of graded cocommutative Hopf algebras $\Lambda \to \Gamma \to \Gamma//\Lambda$, where $\Lambda$ is also commutative. This last commutativity condition is not, in fact, necessary, as forthcoming work of Singer [21] shows.

## 4. Reduction of Theorem 1.2

In this section we state a variant of Theorem 1.2 (namely, Theorem 4.1), and we show that Theorem 1.2 follows from it.

Let $D$ be the sub-Hopf algebra of $A$ generated by the $P_t^s$'s with $s < t$; in other words, $D$ is dual to the following quotient Hopf algebra of $A_*$:

$$D_* = \mathbf{F}_2[\xi_1, \xi_2, \xi_3, \ldots]/(\xi_1^2, \xi_2^4, \xi_3^8, \ldots, \xi_n^{2^n}, \ldots).$$

By Proposition 2.1(c), every elementary sub-Hopf algebra of $A$ sits inside of $D$. By Proposition 2.4(a), $D$ is a normal sub-Hopf algebra of $A$, so there is an action of $A//D$ on $H^*(D; \mathbf{F}_2)$.

We have the following result.

THEOREM 4.1. *The restriction map $H^*(A; \mathbf{F}_2) \to H^*(D; \mathbf{F}_2)$ factors through*

$$\varphi : H^*(A; \mathbf{F}_2) \to H^*(D; \mathbf{F}_2)^{A//D},$$

*and $\varphi$ is an F-isomorphism.*

In this section, we reduce Theorem 1.2 to Theorem 4.1; we prove Theorem 4.1 in Section 5. The first step in the reduction is the following.

PROPOSITION 4.2. *There is an F-isomorphism*

$$f : H^*(D; \mathbf{F}_2) \to \varprojlim_{\mathcal{Q}} H^*(E; \mathbf{F}_2).$$

*Proof.* We need to show two things: every element in the kernel of $f$ is nilpotent, and we can lift some $2^n$-th power of any element in the codomain of $f$. The first of these statements follows from the main theorem of [4]—restriction to the elementary sub-Hopf algebras of $A$ (and hence of $D$) detects nilpotence.

The second statement follows, almost directly, from a result of Hopkins and Smith [6, Theorem 4.13], stated as Theorem 4.3 below, combined with the fact that the Hopf algebra $D$ is an inverse limit of finite-dimensional Hopf algebras. We explain. For each integer $r \geq 1$, we let $B(r)$ be the sub-Hopf algebra of $D$ which is dual to the following quotient Hopf algebra of $D_*$:

$$B(r)_* = D_*/(\xi_1, \xi_2, \ldots, \xi_r).$$

Then each $B(r)$ is normal in $D$, $D//B(r)$ is finite-dimensional, and we have $D = \varprojlim_r D//B(r)$. Note that this inverse limit stabilizes in any given degree.

We claim that given any object $E$ of $\mathcal{Q}$ and given any $y \in H^*(E; \mathbf{F}_2)$, some $2^n$-th power of $y$ is in the image of the restriction map $H^*(D; \mathbf{F}_2) \to H^*(E; \mathbf{F}_2)$. We have a normal sub-Hopf algebra $E(r) = E \cap B(r)$ of $E$, and $E//E(r)$ is a sub-Hopf algebra of $D//B(r)$. We have the following diagram of Hopf algebra extensions, in which the top row is contained, term-by-term, in the bottom:

$$\begin{array}{ccccc} E(r) & \to & E & \to & E//E(r) \\ \downarrow & & \downarrow & & \downarrow \\ B(r) & \to & D & \to & D//B(r). \end{array}$$

Furthermore, $E = \varprojlim E//E(r)$, so any $y \in H^*(E; \mathbf{F}_2)$ is represented by a class $\overline{y} \in H^*(E//E(r); \mathbf{F}_2)$, for some $r$.

Since $E//E(r) \subseteq D//B(r)$ is an inclusion of finite-dimensional graded connected cocommutative Hopf algebras, then by Theorem 4.3, there is an integer $N$ so that the class $\overline{y}^{2^N}$ lifts to $H^*(D//B(r); \mathbf{F}_2)$. Hence $y^{2^N}$ lifts to $H^*(D; \mathbf{F}_2)$. □

We have just used the following theorem:

THEOREM 4.3 (Theorem 4.13 in [6]). *Suppose that $\Lambda \subseteq \Gamma$ is an inclusion of finite-dimensional graded connected cocommutative Hopf algebras over a field $k$ of characteristic $p > 0$. For any $\lambda \in H^*(\Lambda; k)$, there is a number $N$ so that $\lambda^{p^N}$ is in the image of the restriction map $H^*(\Gamma; k) \to H^*(\Lambda; k)$.*

*Proof of Theorem* 1.2, *assuming Theorem* 4.1. By Proposition 4.2, there is an $F$-isomorphism

$$f : H^*(D; \mathbf{F}_2) \to \varprojlim_{\mathcal{Q}} H^*(E; \mathbf{F}_2) = \varprojlim_{\mathcal{Q}_0} H^*(E; \mathbf{F}_2).$$

We need to show that this induces an $F$-isomorphism on the invariants; then composing with the $F$-isomorphism in Theorem 4.1 gives the desired result (it is clear that the composite of two $F$-isomorphisms is again an $F$-isomorphism).

Since each restriction map $H^*(D; \mathbf{F}_2) \to H^*(E; \mathbf{F}_2)$ is an $A$-algebra map when $E$ is normal, so is the assembly of these into the map $f$. So we have an induced map on the invariants:

$$\tilde{f} : (H^*(D; \mathbf{F}_2))^A \to \left( \varprojlim H^*(E; \mathbf{F}_2) \right)^A.$$

Theorem 4.1 says that the restriction map $H^*(A; \mathbf{F}_2) \to H^*(D; \mathbf{F}_2)$ factors through the invariants, and composing $\tilde{f}$ with this yields the factorization of $\rho : H^*(A; \mathbf{F}_2) \to \varprojlim H^*(E; \mathbf{F}_2)$ through the invariants.

Any $x$ in the kernel of $\tilde{f}$ is also in the kernel of $f$, and hence is nilpotent. Given any $y \in (\varprojlim_{\mathcal{Q}} H^*(E; \mathbf{F}_2))^A$, we know that $y^{2^n}$ is in $\text{im}(f)$ for some $n$; by Lemma 4.4 below, we may conclude that $y^{2^{n+N}} \in \text{im}(\tilde{f})$ for some $N$. (To apply the lemma, we need to know that for every $z \in H^*(D; \mathbf{F}_2)$, $\text{Sq}^{2^k}(z)$ is zero for all but finitely many values of $k$. The action of $A$ decreases degrees, by Proposition 2.3, so this is clear.)

This finishes the proof of Theorem 1.2, given Theorem 4.1. □

We have used the following lemma. We say that an $A$-module $M$ is *locally finite* (also known as "tame") if for every element $m \in M$, the submodule generated by $m$ is finite-dimensional.

LEMMA 4.4. *Suppose that $R$ and $S$ are commutative $A$-algebras, with an $A$-algebra map $f : R \to S$ that detects nilpotence: every $x \in \ker f$ is nilpotent. Suppose also that $R$ is locally finite as an $A$-module. Given $z \in R$ so that $f(z) \in S$ is invariant under the $A$-action, then $z^{2^N}$ is also invariant, for some $N$.*

*Proof.* Since $R$ is locally finite, then $\text{Sq}^n(z) = 0$ for all but finitely many values of $n$. Since $f$ is an $A$-module map and $f(z)$ is invariant, then we also know that $f(\text{Sq}^n(z)) = 0$ for all $n$; therefore, each $\text{Sq}^n(z)$ is nilpotent. Recall that the total Steenrod square $\text{Sq} = \sum_{i \geq 0} \text{Sq}^i$ is an algebra map, and to say that an element $y$ is invariant under the $A$-action is to say that $\text{Sq}(y) = y$. In our case, $\text{Sq}(z)$ has only finitely many nonzero terms, and all but the first ($\text{Sq}^0 z = z$) are nilpotent. Hence, by commutativity, $\text{Sq}(z) - z$ is nilpotent; i.e., $z^{2^N}$ is invariant for some $N$. □

## 5. Proof of Theorem 4.1

In this section we prove Theorem 4.1. First we show that the restriction map $H^*(A; \mathbf{F}_2) \to H^*(D; \mathbf{F}_2)$ factors through the invariants $H^*(D; \mathbf{F}_2)^A$, and then we show that the resulting map is an $F$-isomorphism.

$D$ is a normal sub-Hopf algebra of $A$, so by Proposition 2.3(b), there is a spectral sequence

$$E_2^{s,t,u} = H^{s,u}(A//D; H^{t,*}(D; \mathbf{F}_2)) \Rightarrow H^{s+t,u}(A; \mathbf{F}_2).$$

It is well-known, and easy to verify from any reasonable construction of the spectral sequence, that the edge homomorphism

$$H^*(A; \mathbf{F}_2) \to E_2^{0,*} = H^*(D; \mathbf{F}_2)^{A//D}$$

factors the restriction map (see [1, p. 114], for instance). This gives us the desired map

$$\varphi : H^*(A; \mathbf{F}_2) \to H^*(D; \mathbf{F}_2)^{A//D} = H^*(D; \mathbf{F}_2)^A.$$

We have to show two things: that every element of $\ker \varphi$ is nilpotent, and that for every $y \in H^*(D; \mathbf{F}_2)^A$, there is an integer $n$ so that $y^{2^n} \in \text{im} \, \varphi$. (The first of these follows from [15, Theorem 3.1]; see also [4, Theorem 1.1]. We provide a proof here, anyway.) These two verifications have many common features, so first we establish some tools that we will use in both parts.

The idea of the proof is to show that there is a nice vanishing plane at some $E_r$-term of the spectral sequence associated to the extension $D \to A \to A//D$. It turns out to be easier, though, to introduce some intermediate Hopf algebras between $A$ and $D$, and find vanishing planes in the spectral sequences for each of the intermediate extensions. Once the vanishing planes have been established, the rest of the proof is easy.

For $n \geq 1$ we define sub-Hopf algebras $D(n)$ of $A$ by

$$D(n) = \left(\mathbf{F}_2[\xi_1, \xi_2, \xi_3, \ldots]/(\xi_1^2, \xi_2^4, \ldots, \xi_n^{2^n})\right)^*.$$

We let $D(0) = A$. This gives us a diagram of Hopf algebra inclusions:

$$\ldots \hookrightarrow D(2) \hookrightarrow D(1) \hookrightarrow D(0) = A.$$

$D$ is the inverse limit of this diagram. Taking cohomology gives

$$H^*(D(0); \mathbf{F}_2) \to H^*(D(1); \mathbf{F}_2) \to H^*(D(2); \mathbf{F}_2) \to \ldots .$$

By Proposition 2.4, each $D(n)$ is a normal sub-Hopf algebra of $A$; hence there is an action of $A$ on $H^*(D(n); \mathbf{F}_2)$.

LEMMA 5.1.

(a) *There is an isomorphism* $H^*(D; \mathbf{F}_2) \cong \varinjlim H^*(D(n); \mathbf{F}_2)$.

(b) *This induces an isomorphism* $H^*(D; \mathbf{F}_2)^A \cong \varinjlim (H^*(D(n); \mathbf{F}_2))^A$.

*Proof.* Since $D$ is the inverse limit of the $D(n)$'s, and since the diagrams

$$\ldots \hookrightarrow D(3) \hookrightarrow D(2) \hookrightarrow D(1) \hookrightarrow D(0),$$
$$H^*(D(0); \mathbf{F}_2) \to H^*(D(1); \mathbf{F}_2) \to H^*(D(2); \mathbf{F}_2) \to \ldots,$$

stabilize (and are finite-dimensional) in each degree, then part (a) is clear.

Part (b): Since the action of $\mathrm{Sq}^k \in A$ takes $H^{s,t}(B; \mathbf{F}_2)$ to $H^{s,t-k}(B; \mathbf{F}_2)$ (for $B = D$ or $B = D(n)$), then the action of $A$ on any element of $H^*(D; \mathbf{F}_2) \cong \varinjlim H^*(D(n); \mathbf{F}_2)$ is completely determined at some finite stage of the direct system. Hence the ring of invariants of the direct limit is equal to the direct limit of the invariants. □

This lemma lets us translate information about $H^*(D; \mathbf{F}_2)$ to $H^*(D(n); \mathbf{F}_2)$ for some $n \gg 0$. We also want to know how to go from $D(n)$ to $D(n-1)$ (and then apply downward induction on $n$, ending at $D(0) = A$). There is a Hopf algebra extension

(**) $$D(n) \to D(n-1) \to D(n-1)//D(n),$$

and the quotient Hopf algebra is easy to describe: it is dual to $\mathbf{F}_2[\xi_n^{2^n}]$ with $\xi_n^{2^n}$ primitive, so $D(n-1)//D(n)$ is isomorphic to an exterior algebra:

$$D(n-1)//D(n) \cong E(P_n^n, P_n^{n+1}, P_n^{n+2}, \ldots).$$

(Recall from Section 2 that $P_n^{n+s}$ is the element of $A$ dual to $\xi_n^{2^{n+s}}$.) Anyway, by Proposition 2.3, for any $D(n-1)$-module $M$ there is a spectral sequence

$$E_2^{s,t,u}(M) = H^{s,u}(D(n-1)//D(n); H^{t,*}(D(n); M)) \Rightarrow H^{s+t,u}(D(n-1); M).$$

Our main tool for the completion of the proof is the following.

LEMMA 5.2. *Fix an integer $N$. For some $r$ and some $c$, there is a vanishing plane "of slope $N$" at the $E_r$-term of the spectral sequence: $E_r^{s,t,u}(\mathbf{F}_2) = 0$ when $Ns + t - u > c$.*

To prove this lemma, we show that for a certain finite module $M$, we have the appropriate vanishing plane at $E_2^{***}(M)$; then we use Theorem 3.2 and Proposition 3.6 to induce a vanishing plane (with the same normal vector, but possibly different intercept) at $E_r^{***}(\mathbf{F}_2)$ for some $r$.

LEMMA 5.3. *Fix an integer $j > n$, and let $M_j$ be the $D(n-1)//D(n)$-module dual to $E[P_n^n, P_n^{n+1}, \ldots, P_n^{j-1}]$. Then $E_2^{s,t,u}(M_j) = 0$ when*

$$2^j(2^n - 1)s + t - u > 0.$$

(So given an integer $N$, if we choose $j$ so that $2^j(2^n - 1) > N$, then we will have $E_2^{s,t,u}(M_j) = 0$ for $Ns + t - u > 0$.)

*Proof.* $M_j$ is a trivial $D(n)$-module (by definition), and as algebras, we have

$$D(n-1)//D(n) \cong M_j^* \otimes E[P_n^j, P_n^{j+1}, P_n^{j+2}, \ldots].$$

Hence the $E_2$-term of the spectral sequence for $M_j$ looks like

$$\begin{aligned} E_2^{s,t,u}(M_j) &\cong H^{s,u}(D(n-1)//D(n); H^{t,*}(D(n); M_j)) \\ &\cong H^{s,u}(D(n-1)//D(n); M_j \otimes H^{t,*}(D(n); \mathbf{F}_2)) \\ &\cong H^{s,u}(E[P_n^j, P_n^{j+1}, \ldots]; H^{t,*}(D(n); \mathbf{F}_2)). \end{aligned}$$

The Hopf algebra $E[P_n^j, P_n^{j+1}, \ldots]$ is $2^j(2^n - 1)$-connected, so if $L$ is a module which is zero below degree $t$, then

$$H^{s,u}(E[P_n^j, P_n^{j+1}, \ldots]; L) = 0$$

when $u < 2^j(2^n - 1)s + t$. Now note that $H^{t,*}(D(n); \mathbf{F}_2)$ is zero below degree $t$. □

We want to apply Theorem 3.2 to get a vanishing plane at some term of the spectral sequence $E_r^{***}(\mathbf{F}_2)$; first we describe how to build $M_j$ out of $\mathbf{F}_2$, inductively.

LEMMA 5.4. *For $j > n$, let $M_j$ be as in Lemma 5.3. Let $M_n = \mathbf{F}_2$. For each $j \geq n$, there is a short exact sequence of $D(n-1)//D(n)$-modules*

$$0 \to \Sigma^{2^j(2^n-1)} M_j \to M_{j+1} \to M_j \to 0.$$

*This represents the element $h_{nj} \otimes 1_{M_j} \in \operatorname{Ext}^{**}_{D(n-1)}(M_j, M_j)$, where*

$$h_{nj} = [\xi_n^{2^j}] \in H^*(D(n-1); \mathbf{F}_2).$$

*Proof.* This is clear. □

*Proof of Lemma 5.2.* Lin has shown in [8, 3.2] that when $j \geq n$, then the element $h_{nj} \in H^*(D(n-1); \mathbf{F}_2)$ is nilpotent. Therefore Lemma 5.2 follows from the above lemmas, together with Theorem 3.2 and Proposition 3.6. □

By the way, Lemma 5.2 has the following obvious corollary; we do not use this in this paper, but it might be useful elsewhere.

COROLLARY 5.5. *Fix an integer $N$ and a finite-dimensional $A$-module $M$, and consider the Lyndon-Hochschild-Serre spectral sequence $E_*^{s,t,u}(M)$ associated to the extension (\*\*). For some $r$ and some $c$, there is a vanishing plane at the $E_r$-term of this spectral sequence: $E_r^{s,t,u}(M) = 0$ when $Ns + t - u > c$.*

*Proof.* This follows from Theorem 3.2 and induction on $\dim M$. □

Now we finish the proof of Theorem 4.1. We have two conditions to verify, and we devote a lemma to each one.

LEMMA 5.6. *Every $x \in \ker(\varphi)$ is nilpotent.*

*Proof.* Fix $x \in H^*(A; \mathbf{F}_2)$ so that $\varphi(x) = 0$; we want to show that $x$ is nilpotent. As mentioned above, this has been proved elsewhere; see [15, Theorem 3.1] and [4, Theorem 1.1]. We present a proof here anyway, because we have done most of the work already in setting up the other half of the proof of Theorem 4.1.

Given Hopf algebras $C \subseteq B$ over a field $k$, we write $\operatorname{res}_{B,C} : H^*(B;k) \to H^*(C;k)$ for the restriction map.

We have an element $x \in H^*(A; \mathbf{F}_2)$ so that $\varphi(x) = 0$. Then $\operatorname{res}_{A,D}(x) = 0$, so by Lemma 5.1(a), $\operatorname{res}_{A,D(n)}(x) = 0$ for some $n$. Since $D(0) = A$ and $\operatorname{res}_{A,D(0)}$ is the identity map on $H^*(A; \mathbf{F}_2)$, it suffices to show that if for some $n$, $\operatorname{res}_{A,D(n)}(x) = 0$, then $\operatorname{res}_{A,D(n-1)}(x^j) = 0$ for some $j$.

Consider the spectral sequence

$$E_2^{s,t,u}(\mathbf{F}_2) = H^{s,u}(D(n-1)//D(n); H^{t,*}(D(n); \mathbf{F}_2)) \Rightarrow H^{s+t,u}(D(n-1); \mathbf{F}_2)$$

associated to the extension

$$D(n) \to D(n-1) \to D(n-1)//D(n).$$

Write $z$ for $\operatorname{res}_{A,D(n-1)}(x)$, so that $\operatorname{res}_{D(n-1),D(n)}(z) = \operatorname{res}_{A,D(n)}(x) = 0$. Because of this, $z$ must be represented by a class $\tilde{z} \in E_2^{p,q,v}$ with $p > 0$.

Choose an integer $N$ so that $Np + q - v > 0$. By Lemma 5.2, for some $r$ and $c$ we have $E_r^{s,t,u}(\mathbf{F}_2) = 0$ when $Ns + t - u > c$. Since this is a spectral sequence of commutative algebras, we can find a $j$ so that $\tilde{z}^j$ survives to $E_r$, and also so that
$$Npj + qj - vj > c.$$

Hence at this $E_r$-term, $\tilde{z}^j$ must be zero. This means that $\tilde{z}^j = 0$ at $E_\infty$, so that $\tilde{z}^j = 0$ in the abutment, modulo higher filtrations. But the vanishing plane at $E_r$ also exists at $E_\infty$, and one can see that the higher filtrations are all above the vanishing plane, and hence zero. So $\tilde{z}^j = 0$ in the abutment; hence $z^j = 0$ in $H^*(D(n-1); \mathbf{F}_2)$, which is what we wanted to show. □

LEMMA 5.7. *Every $y \in \operatorname{coker}(\varphi)$ is nilpotent.*

*Proof.* Fix $y \in H^*(D; \mathbf{F}_2)^{A//D}$. We want to show that there is an integer $m$ so that $y^{2^m} \in \operatorname{im} \varphi$. By Lemma 5.1(b), there is an $n$ so that $y$ lifts to $H^*(D(n); \mathbf{F}_2)^{A//D(n)}$. Now we show that some power of $y$ lifts to
$$H^*(D(n-1); \mathbf{F}_2)^{A//D(n-1)}.$$

Since $D(0) = A$, then downward induction on $n$ will finish the proof.

Since $y \in H^*(D(n); \mathbf{F}_2)$ is invariant under the $A//D(n)$-action, then it is also invariant under the action of $D(n-1)//D(n)$ (since the latter quotient is a sub-Hopf algebra of the former). So $y$ represents a class at the $E_2$-term of the spectral sequence
$$E_2^{s,t,u}(\mathbf{F}_2) = H^{s,u}(D(n-1)//D(n); H^{t,*}(D(n); \mathbf{F}_2)) \Rightarrow H^{s+t,u}(D(n-1); \mathbf{F}_2).$$

Also, $y$ lies in the $(t,u)$-plane, say $y \in E_2^{0,q,v}$. Choose $N$ large enough so that $N + q - v - 1$ and $N - 1$ are both positive. By Lemma 5.2, we know that for some $c$ and $r$, we have $E_r^{s,t,u} = 0$ for $Ns + t - u > c$.

Choose $i$ so that
$$2^i > \max(r - 2, \frac{c - N}{N + q - v - 1}).$$

We claim that $y^{2^i}$ is a permanent cycle. By Lemma 3.7, the possible differentials on $y^{2^i}$ are
$$d_{j+1}(y^{2^i}) \in E_{j+1}^{j+1, 2^i q - j, 2^i v},$$

for $j \geq 2^i$, so we show that the targets of these differentials lie above the vanishing plane. (Since $2^i$ was chosen to be at least $r - 1$, then the vanishing plane is present at the $E_{j+1}$-term for all $j \geq 2^i$.) We just have to verify the inequality specified by the vanishing plane:

$$\begin{aligned} N(j+1) + 2^i q - j - 2^i v &= (N-1)j + N + 2^i(q-v) \\ &\geq (N-1)2^i + N + 2^i(q-v) \\ &= 2^i(N+q-v-1) + N \\ &> \frac{c-N}{N+q-v-1}(N+q-v-1) + N \\ &= c. \end{aligned}$$

Hence all of the differentials on $y^{2^i}$ vanish, so it is a permanent cycle. For degree reasons, it cannot be a boundary; hence it gives a nonzero element of $E_\infty$, and hence a nonzero element of $H^*(D(n-1); \mathbf{F}_2)$. Lastly, Lemma 4.4 tells us that some $2^k$-th power of $y^{2^i}$ is invariant in $H^*(D(n-1); \mathbf{F}_2)$. This completes the inductive step, hence the proof of Lemma 5.7, and hence the proof of Theorem 4.1. □

MASSACHUSETTS INSTITUTE OF TECHNOLOGY, CAMBRIDGE, MA
*Current address*: UNIVERSITY OF NOTRE DAME, NOTRE DAME, IN
*E-mail address*: palmieri@member.ams.org